\documentclass[12pt,draftcls,onecolumn]{IEEEtran}
\IEEEoverridecommandlockouts

\usepackage{amsmath}
\usepackage{graphicx}
\usepackage{cite}
\usepackage{authblk}
\usepackage{amssymb}
\usepackage{subfigure}
\usepackage{multirow}
\usepackage{color}
\usepackage{tikz} \usetikzlibrary{arrows,calc,decorations.pathreplacing}
\usepackage{soul} \setstcolor{red}
\usepackage{psfrag}
         
\newtheorem{theorem}{Theorem}

\title{Efficient Communication Leader-Wingman Flight Formation Control}

\author{Eloy Garcia\thanks{Eloy Garcia is with the Control Science Center of Excellence, Air Force Research Laboratory, Wright-Patterson AFB, OH 45433. \ttfamily{eloy.garcia.2@us.af.mil}}}

\begin{document}
\maketitle 
\begin{abstract}
This paper considers autonomous flight formation of fixed-wing air vehicles. Due to several issues such as limitations on communication bandwidth and low observability requirements, an asynchronous control and communication approach is designed in order to reduce the number of updates between the vehicles. With the implementation of an event-triggered controller, it is possible to obtain a more efficient communication scheme while also achieving the convergence properties of the desired formation. 
\end{abstract}

\section{Introduction} \label{sec:intro}
Autonomous flight formation represents a specific task within autonomous control of aircraft or groups of air vehicles \cite{Pachter98}. The objective of this problem is for an Unmanned Aerial Vehicle (UAV) labeled as the wingman to autonomously achieve a commanded formation relative to the manned or unmanned  leader aircraft. One of the first reported studies on autonomous leader-wingman flight formation appeared in  \cite{Pachter94}. Different approaches have been subsequently developed by different authors \cite{Boskovic02,Innocenti03,Proud99}.

The authors of \cite{Giulietti04} analyzed the formation flight control problem of one leader and several wingmen and provided simulations using nonlinear models of aircraft. Controller design and experimental results concerning leader-wingman formation flight were presented in \cite{Seanor04}.

The controller design offered in the current paper differs from \cite{Seanor04} mainly by addressing the low communication rate between aircraft. We consider an event-triggered control and communication approach 
\cite{AntaTabuada10,Aragues12,Garcia12,Eqtami11IFAC,DePersis13,Demir12a,Mazo11TAC,yu2013model,Guinaldo11,Cervin08CDC,Garcia12Parameter,Guinaldo12CDC,LiuHill12,Garcia14CDCetof,WangLemmon08a,Mccourt14,Adaldo15,nowzari19,Tallapragada12IFAC} in order to reduce communication rate between aircraft.

Event-triggered control and stabilization of dynamical systems was considered first in \cite{Astrom02,Astrom08} and later formalized by Tabuada and Wang in \cite{Wang06,Tabuada07}. 
The model-based approach for control of networked systems and systems with limited feedback has been addressed in \cite{Garcia13scl,Garcia13MED,Garcia12MED,GarciaBook} where periodic communication was used. 
The model-based event-triggered (MB-ET) control combines the previous approaches \cite{lehmann2012event,Garcia20arxiv} in order to maximize the time intervals for which a control system is able to operate without feedback information. In this paper we take advantage of the MB-ET control architecture in order stabilize the leader-wingman formation. Different from \cite{Seanor04,Pachter98}, where the formation control algorithm requires continuous communication between leader and wingman, the controller in this paper reduces frequency of communication updates. In particular, feedback updates are transmitted in an asynchronous manner and only when it is necessary to do so. This requirement is evaluated by means of an event-triggered function which evaluates the size of the state error relative to the current state of the system.

In the following we state the formation control problem in Section \ref{sec:PS}. The nonlinear dynamics of the leader-wingman formation control problem are linearized in Section \ref{sec:LC} and a PI controller is proposed as well.
In Section \ref{sec:example}, the MB-ET control architecture is employed and an illustrative example is presented. Conclusions follow in Section \ref{sec:concl}.

\begin{figure}[htb]\centering
\begin{tikzpicture}[>=stealth]
	\tikzset{
	   missile/.pic = {  
	      \draw [pic actions] (0,0)
	         -- (-4mm,3.5mm)
	         -- ++(8.4mm,-3.5mm)  
	         -- ++(-8.4mm,-3.5mm)
	         -- (0,0);
	   }
	}
	\coordinate (D) at (1,2);  
	\coordinate (A) at (6,4);
	\draw [->] (D) -- ++ (0:1.5)
			(D) ++(0:.6) arc [start angle=0, end angle=50, radius=.6cm]
			(D) ++(0:.9) arc [start angle=0, end angle=22, radius=.9cm]
		        (1.85,2.55) node{$\psi_W$}
			(2.3,3.25) node{$V_W$}
			(2.15,2.25) node{$\lambda$}
			(3.5,3.25) node{$R$}
			(1.8,1.05) node{$y$}
			(4.55,2.1) node{$x$}
		        (D) -- ++(50:1.5);
	\draw (D) pic [rotate=50,fill=cyan,scale=.7] {missile};  
	\draw [->] (D) -- ++(21.8:5);
	\draw [->] (D) -- ++(-40:2.482)
	   	 -- ++(50:4.4);
	 \draw [->] (D) -- ++(-40:2.482);
	\draw (D) node [below left] {$Wing$};
	\draw [->] (A) -- ++ (0:1.5)
			(A) ++(0:.7) arc [start angle=0, end angle=70, radius=.7cm]
						(6.8,4.65) node{$\psi_L$}
						(6.8,5.35) node{$V_L$}
		        (A) -- ++(70:1.5);
	\draw (A) pic [rotate=70,fill=brown,scale=.7] {missile};
	\draw (A) node [below right] {$Lead$};
	\node [right] (X) at (0,5) {$X$};
	\draw [->](0,-.4) -- (0,5);
	\node [above] (Y) at (7,0) {$Y$};
	\draw [->](-1,0) -- (7,0);
 \end{tikzpicture}
\caption{Leader-wingman flight formation}
\label{fig:problem description}
\end{figure}
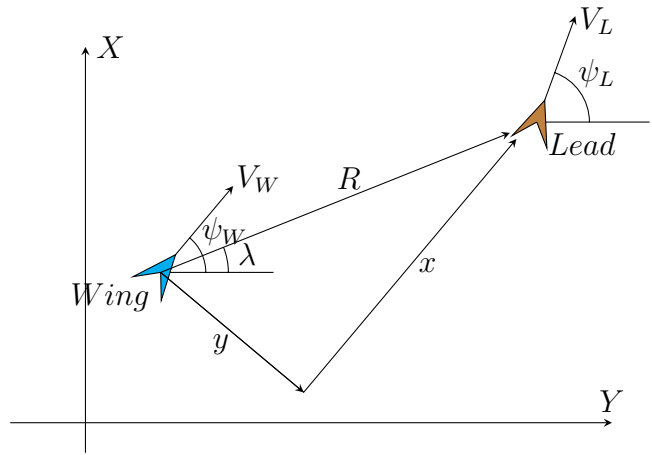


\section{Problem statement}   \label{sec:PS}
We consider the following continuous-time dynamics for the leader
\begin{align}
 \left.
	\begin{array}{l l}
      \dot{x}_L=V_L \cos\psi_L,   \ \ \ \ \ \ \ \     \dot{V}_L=\frac{1}{\tau_{V_L}}(V_{L_c}-V_L)   \\
        \dot{y}_L=V_L \sin\psi_L,   \ \ \ \ \ \ \ \     \dot{\psi}_L=\frac{1}{\tau_{\psi_L}}(\psi_{L_c}-\psi_L).  
	\end{array}  \right.   \label{eq:DynamicsL}
\end{align}
Similarly, the wingman dynamics are given by
\begin{align}
 \left.
	\begin{array}{l l}
      \dot{x}_W=V_W \cos\psi_W,   \ \ \ \ \ \ \ \     \dot{V}_W=\frac{1}{\tau_{V_W}}(V_{W_c}-V_W)   \\
        \dot{y}_W=V_W \sin\psi_W,   \ \ \ \ \ \ \ \     \dot{\psi}_W=\frac{1}{\tau_{\psi_W}}(\psi_{W_c}-\psi_W).  
	\end{array}  \right.   \label{eq:DynamicsW}
\end{align}
where $\tau_i>0$, for $i=\tau_{V_L},\tau_{V_W},\tau_{\psi_L},\tau_{\psi_W}$, are the corresponding time constants. It is assumed that appropriate autopilots are implemented in order to obtain the decoupled $V-\psi$ dynamics \eqref{eq:DynamicsL}-\eqref{eq:DynamicsW}  where both aircraft are controlled by means of their corresponding speed and heading commands. The leader commands $V_{L_c}$ and $\psi_{L_c}$ are considered as external commands which are preset or given by a human pilot. The control inputs for the autonomous wingman are $V_{W_c}$ and $\psi_{W_c}$. 

The problem is to design the feedback controllers in order for the wingman to achieve a desired formation with respect to the leader. A desired formation is specified by $x_r$ and $y_r$ which represent longitudinal and lateral separations with respect to the leader current position and along the current wingman's velocity vector, see Fig. \ref{fig:problem description}. The longitudinal separation is given by $x=R\cos(\psi_W-\lambda)$ and the lateral separation is given by $x=R\sin(\psi_W-\lambda)$, where $R=\sqrt{(x_L-x_W)^2+(y_L-y_W)^2}$ and $\lambda=\arctan(\frac{y_L-y_W}{x_L-x_W})$. The rate of change of $R$ and $\lambda$ are given by the following expressions
\begin{align}
 \left.
	\begin{array}{l l}
      \dot{R}= V_L\cos(\psi_L-\lambda) -V_W\cos(\psi_W-\lambda)  \\
        \dot{\lambda}= \frac{1}{R} [V_L\sin(\psi_L-\lambda) -V_W\sin(\psi_W-\lambda)].  
	\end{array}  \right.   \label{eq:RLambda}
\end{align}
We also obtain the rate of change of the separations $x$ and $y$ and they are given by 
\begin{align}
 \left.
	\begin{array}{l l}
      \dot{x}= V_L\cos\psi_e -V_W - y \dot{\psi}_W  \\
        \dot{y}=  x \dot{\psi}_W - V_L\sin\psi_e
	\end{array}  \right.   \label{eq:DynamicsXY}
\end{align}
where $\psi_e=\psi_L-\psi_W$ is the heading error. 

\section{Linearization}   \label{sec:LC}

In this section the system dynamics are linearized around $\psi_e\approx 0$ and $V_L\approx V_W \approx V_n$, where $V_n$ is the nominal speed.
 The linearized error dynamics are given by
\begin{align}
 \left.
	\begin{array}{l l}
      \dot{x}=V_L  -V_W  +\frac{y_r}{\tau_{\psi_W}} \psi_W -\frac{y_r}{\tau_{\psi_W}} \psi_{W_c}    \\
  \dot{V}_W =   -\frac{1}{\tau_{V_W}} V_W +\frac{1}{\tau_{V_W}} V_{W_c}  \\
  \dot{V}_L =   -\frac{1}{\tau_{V_L}} V_L +\frac{1}{\tau_{V_L}} V_{L_c}  \\
   \dot{y}=  - V_n \psi_L   + \big(V_n - \frac{x_r}{\tau_{\psi_W}}  \big) \psi_W +  \frac{x_r}{\tau_{\psi_W}} \psi_{W_c}  \\
   \dot{\psi}_W =   -\frac{1}{\tau_{\psi_W}} \psi_W +\frac{1}{\tau_{\psi_W}} \psi_{W_c}  \\
  \dot{\psi}_L =   -\frac{1}{\tau_{\psi_L}} \psi_L +\frac{1}{\tau_{\psi_L}} \psi_{L_c}  
	\end{array}  \right.   \label{eq:DynamicsLinearCT}
\end{align}
where the error variables are: $x:=\Delta x =x-x_r$, $y:=\Delta y =y-y_r$, $V_L:=\Delta V_L= V_L-V_n$, and $V_W:=\Delta V_W= V_W-V_n$. 

The problem to be solved in this paper is to determine the control commands $V_{W_c}$ and $\psi_{W_c}$ such that the wingman achieves a formation relative to the leader specified by the parameters $x_r$ and $y_r$. In addition, the formation has to be achieved and maintained in the presence of intermittent communication and using an event-triggered control and communication approach.

Define the state $\textbf{x}=[x \ V_W \ V_L \ y \ \psi_W \ \psi_L ]^T$, the control input $\textbf{u}=[V_{W_c} \ \psi_{W_c}  ]^T$, and the exogenous input $\textbf{d}=[V_{L_c} \ \psi_{L_c}  ]^T$. Then, the linearized dynamics in compact form are given by
\begin{align}
  \dot{\textbf{x}}(t)&=A\textbf{x}(t) +B \textbf{u}(t) +E\textbf{d}(t)   \label{eq:compact}
  \end{align}
  where  
 \begin{align}
  A&= \begin{bmatrix}
	A_{11}&A_{12}  \\
	0&A_{22} 
\end{bmatrix} \\ \qquad
   B&= \begin{bmatrix}
	B_1  \\
	B_2 
\end{bmatrix} \\
 \qquad
E&= \begin{bmatrix}
	E_1  \\
	E_2 
\end{bmatrix}  \label{eq:Lmatrices}
\end{align}
and each block of matrix $A$ is given by
\begin{align}
  A_{11}&= \begin{bmatrix}
	0 &  -1  &   1    \\
	0 & -\frac{1}{\tau_{V_W}}  & 0\\
	 0 & 0  &   -\frac{1}{\tau_{V_L}}
\end{bmatrix}  \nonumber \\ 
A_{12}&= \begin{bmatrix}
	0 &  \frac{y_r}{\tau_{\psi_W}}  &  0  \\
	0 & 0 & 0  \\
	 0 & 0  &   0
\end{bmatrix}  \nonumber \\ 
A_{22}&= \begin{bmatrix}
	0 & V_n -\frac{x_r}{\tau_{\psi_W}} & -V_n   \\
	0 & -\frac{1}{\tau_{\psi_W}}  & 0\\
	 0 & 0  &    -\frac{1}{\tau_{\psi_L}} 
\end{bmatrix}   \nonumber   
\end{align}
Similarly, the block components of the B matrix are given by
\begin{align}
B_1&= \begin{bmatrix}
	0 &  -\frac{y_r}{\tau_{\psi_W}}    \\
	\frac{1}{\tau_{V_W}} & 0 \\
	 0 &  0
\end{bmatrix}  \nonumber \\
B_2&= \begin{bmatrix}
	0 &  \frac{x_r}{\tau_{\psi_W}}    \\
	0 & \frac{1}{\tau_{\psi_W}}  \\
	 0 &  0
\end{bmatrix}  \nonumber 
\end{align}
Finally, matrix $E$ is given by
\begin{align}
E_1&= \begin{bmatrix}
	0 & 0    \\
	0 & 0 \\
	 \frac{1}{\tau_{V_L}}  &  0
\end{bmatrix}  \nonumber \\
E_2&= \begin{bmatrix}
	0 & 0    \\
	0 & 0 \\
	0 &  \frac{1}{\tau_{\psi_L}}  
\end{bmatrix}  \nonumber 
\end{align}
In reference \cite{Pachter94} a Proportional-Integral (PI) controller for each channel ($x$ and $y$) was proposed.
The PI controller for the $x$-channel is given by
\begin{align}
 \left.
	\begin{array}{l l}
    V_{W_c} = K_{xp}e_x + K_{xI} \int_0^t e_x dt
	\end{array}  \right.   \label{eq:CTPIx}
\end{align}
 where $e_x=K_x x + K_V(V_L-V_W)$. The PI controller for the $y$-channel is given by
\begin{align}
 \left.
	\begin{array}{l l}
    \psi_{W_c} = K_{yp}e_y + K_{yI} \int_0^t e_y dt
	\end{array}  \right.   \label{eq:CTPIy}
\end{align}
 where $e_y=K_y y + K_\psi(\psi_L-\psi_W)$. However, conditions on the control gains to guarantee stability of the formation were not presented in \cite{Pachter94}.
In the remaining of this section, stability conditions are obtained by augmenting the state of the system and analyzing the the two separate channels. The analysis in this section will also be useful in the following section in order to implement a model-based event-triggered controller.

\begin{theorem} 
The linearized system \eqref{eq:compact} with control inputs given by \eqref{eq:CTPIx} and \eqref{eq:CTPIy}
is stable if the control gains in the $x$-channel satisfy the following conditions
\begin{align}
 \left.
	\begin{array}{l l}
       K_{xp}K_V + 1> 0 \\
         K_xK_{xI} > 0  \\
           K_{xp}K_x + K_V K_{xI} > \frac{\tau_{V_W}K_xK_{xI}}{K_{xp}K_V + 1}  
    	\end{array}  \right.   \label{eq:StabCondsx}
\end{align}
and the control gains in the $y$-channel satisfy the following conditions
\begin{align}
 \left.
	\begin{array}{l l}
       K_{yp}K_\psi - x_rK_yK_{yp} + 1> 0 \\
         K_yK_{yI} < 0  \\
           K_{yI}K_\psi - K_y (V_nK_{yp}+x_rK_{yI})  > - \frac{\tau_{\psi_W}V_nK_yK_{yI}}{K_{yp}K_\psi - x_rK_yK_{yp} + 1}  
    	\end{array}  \right.   \label{eq:StabCondsy}
\end{align}
\end{theorem}

\section{MB-ET control and example}   \label{sec:example}

In order to implement the controller \eqref{eq:CTPIx}-\eqref{eq:CTPIy}, the leader needs to transmit measurements of its own state to the wingman. In this section, the goal is to reduce the number of transmission by implementing the MB-ET control framework. For this problem, both the Leader and the Wingman, implement a model of the Leader dynamics. The state of such model is denoted as $\hat{\textbf{x}}=[\hat{x}_L \ \hat{y}_L \ \hat{V}_L \ \hat{\psi}_L]^T$. The state of each model is updated at the update time instants $t_k$, that is, $\hat{\textbf{x}}(t_k)=\textbf{x}(t_k)$, where $\textbf{x}=[x_L \ y_L \ V_L \ \psi_L]^T$ is the state of the Leader. The state error is defined as $\textbf{e}=\hat{\textbf{x}}-\textbf{x}$. The events are generated when $|| \textbf{e}|| \geq \sigma ||\textbf{x}||$, where $0<\sigma<1$.

\begin{figure}
	\begin{center}
		\includegraphics[width=10.0cm,trim=1.5cm .2cm 1.5cm .5cm]{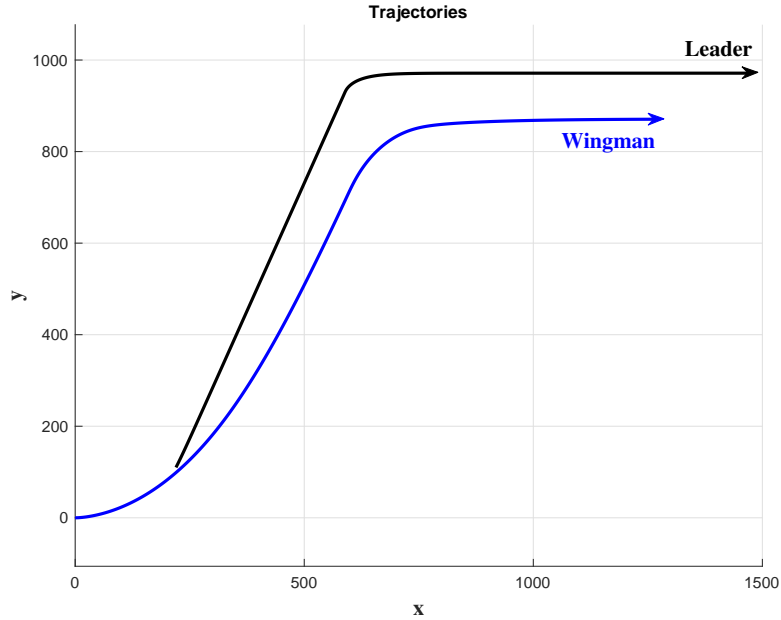}
	\caption{Trajectories of Leader and wingman aircraft.}
	\label{fig:fig1}
	\end{center}
\end{figure}

\begin{figure}
	\begin{center}
		\includegraphics[width=8.4cm,trim=1.5cm .2cm 1.5cm .5cm]{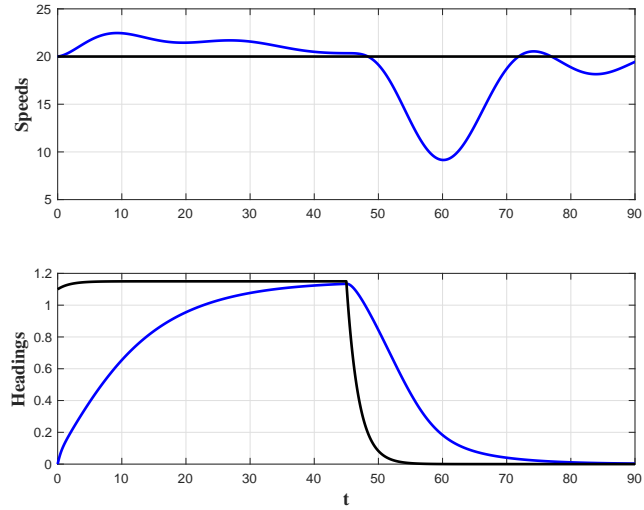}
	\caption{Speeds and headings of aircraft.}
	\label{fig:fig2}
	\end{center}
\end{figure}

\begin{figure}
	\begin{center}
		\includegraphics[width=8.3cm,trim=1.5cm .2cm 1.5cm .5cm]{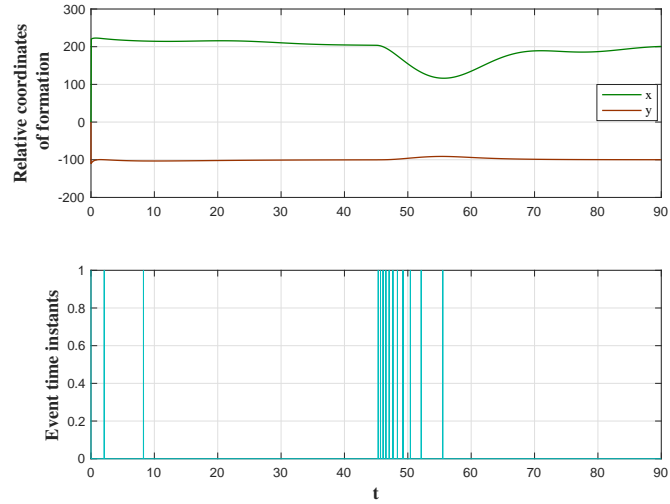}
	\caption{Formation coordinates and event time instants.}
	\label{fig:fig3}
	\end{center}
\end{figure}

\textit{Example}. Consider the Leader's nominal speed $V_n =20$. The desired relative coordinates are $x_r=200$ and $y_r=-100$. The initial conditions of the players are such that the initial errors are $x-x_r=20$ and $y-y_r=-10$.
Initially, the Leader is commanded to implement a heading of $\psi_{L_c}=1.15$ rad. The wingman's initial heading is $\psi_{W_0}=0$ rad. At time $t=45$ sec, the Leader implements a maneuver and changes its heading to $\psi_{L_c}=0$ rad. The resulting trajectories of the Leader and the wingman are shown in Fig. \ref{fig:fig1}. The speeds and headings of the aircraft are shown in Fig. \ref{fig:fig2}. In this example the Leader is not commanded to change its speed but only its heading. The wingman computes its control inputs, its control heading and its own speed, in order to achieve and keep the desired formation. It can be seen in Fig. \ref{fig:fig3} that the desired formation is achieved. The initial error is suppressed and the relative coordinates converge to the desired values $x_r$ and $y_r$.

At the maneuver time $t=45$ sec, the Leader suddenly receives the external command to change it s heading. This maneuver creates a transient error on the formation and the wingman has to act accordingly in order to preserve the formation. The time instants when the Leader transmits a state measurement to the wingman are shown at the bottom of  Fig. \ref{fig:fig3}. It is clear that communication is needed when there are transient formation errors but communication is not as necessary once the aircraft achieve the formation and the external commands remain constant, that is, a new maneuver is not commanded.

\section{Conclusions}   \label{sec:concl}

Communication constraints within the leader-wingman formation control problem were addressed in this paper. The model-based event-triggered (MB-ET) control framework was implemented in order to minimize the number of transmissions from the Leader to the Wingman. The MB-ET control framework combines the model-based control approach with the event-triggered control and communication paradigm which results in significant savings of communication bandwidth.

\bibliographystyle{ieeetran}
\bibliography{ReferencesFF}

\begin{thebibliography}{10}
\providecommand{\url}[1]{#1}
\csname url@samestyle\endcsname
\providecommand{\newblock}{\relax}
\providecommand{\bibinfo}[2]{#2}
\providecommand{\BIBentrySTDinterwordspacing}{\spaceskip=0pt\relax}
\providecommand{\BIBentryALTinterwordstretchfactor}{4}
\providecommand{\BIBentryALTinterwordspacing}{\spaceskip=\fontdimen2\font plus
\BIBentryALTinterwordstretchfactor\fontdimen3\font minus
  \fontdimen4\font\relax}
\providecommand{\BIBforeignlanguage}[2]{{%
\expandafter\ifx\csname l@#1\endcsname\relax
\typeout{** WARNING: IEEEtran.bst: No hyphenation pattern has been}%
\typeout{** loaded for the language `#1'. Using the pattern for}%
\typeout{** the default language instead.}%
\else
\language=\csname l@#1\endcsname
\fi
#2}}
\providecommand{\BIBdecl}{\relax}
\BIBdecl

\bibitem{Pachter98}
M.~Pachter and P.~R. Chandler, ``Challenges of autonomous control,'' \emph{IEEE
  Control Systems}, vol.~18, no.~4, pp. 92--97, 1998.

\bibitem{Pachter94}
M.~Pachter, J.~J. D'Azzo, and J.~L. Dargan, ``Automatic formation flight
  control,'' \emph{AIAA Journal of Guidance, Control, and Dynamics}, vol.~17,
  no.~6, pp. 1380--1383, 1994.

\bibitem{Boskovic02}
J.~D. Boskovic, S.~M. Li, and R.~K. Mehra, ``Formation flight control design in
  the presence of unknown leader commands,'' in \emph{American Control
  Conference}, 2002, pp. 2854--2859.

\bibitem{Innocenti03}
M.~Innocenti, F.~Giulietti, and L.~Pollini, ``Intelligent management control
  for unmanned aircraft navigation and formation,'' 2003, pp. 1--22.

\bibitem{Proud99}
A.~Proud, M.~Pachter, and J.~J. D'Azzo, ``Close formation flight control,'' in
  \emph{AIAA Guidance, Navigation, and Control Conference}, AIAA 99-4207, 1999,
  pp. 1--16.

\bibitem{Giulietti04}
F.~Giulietti and G.~Mengali, ``Dynamics and control of different aircraft
  formation structures,'' \emph{The Aeronautical Journal}, vol. 108, no. 1081,
  pp. 117--124, 2004.

\bibitem{Seanor04}
B.~Seanor, G.~Campa, Y.~Gu, M.~Napolitano, L.~Rowe, and M.~Perhinschi,
  ``Formation flight test results for {UAV} research aircraft models,'' in
  \emph{AIAA 1st Intelligent Systems Technical Conference}, AIAA 2004-6251,
  2010, pp. 1--14.

\bibitem{AntaTabuada10}
A.~Anta and P.~Tabuada, ``To sample or not to sample: Self-triggered control
  for nonlinear systems,'' \emph{IEEE Transactions on Automatic Control},
  vol.~55, no.~9, pp. 2030--2042, 2010.

\bibitem{Aragues12}
R.~Aragues, G.~Shi, D.~V. Dimarogonas, C.~Sagues, and K.~H. Johansson,
  ``Distributed algebraic connectivity estimation for adaptive event-triggered
  consensus,'' in \emph{American Control Conference}, 2012, pp. 32--37.

\bibitem{Garcia12}
E.~Garcia and P.~J. Antsaklis, ``Decentralized model-based event-triggered
  control of networked systems,'' in \emph{American Control Conference}, 2012,
  pp. 6485--6490.

\bibitem{Eqtami11IFAC}
A.~Eqtami, D.~V. Dimarogonas, and K.~J. Kyriakopoulos, ``Event-triggered
  strategies for decentralized model predictive controllers,'' in
  \emph{Proceedings of the IFAC World Congress}, 2011, pp. 10\,068--10\,073.

\bibitem{DePersis13}
C.~D. Persis, R.~Sailer, and F.~Wirth, ``Parsimonious event-triggered
  distributed control: a zeno free approach,'' \emph{Automatica}, vol.~49,
  no.~7, pp. 2116--2124, 2013.

\bibitem{Demir12a}
O.~Demir and J.~Lunze, ``Cooperative control of multi-agent systems with
  event-based communication,'' in \emph{American Control Conference}, 2012, pp.
  4504--4509.

\bibitem{Mazo11TAC}
M.~Mazo and P.~Tabuada, ``Decentralized event-triggered control over wireless
  sensor/actuator networks,'' \emph{IEEE Transactions on Automatic Control},
  vol.~56, no.~10, pp. 2456--2461, 2011.

\bibitem{yu2013model}
H.~Yu, E.~Garcia, and P.~J. Antsaklis, ``Model-based scheduling for networked
  control systems,'' in \emph{2013 American Control Conference}.\hskip 1em plus
  0.5em minus 0.4em\relax IEEE, 2013, pp. 2350--2355.

\bibitem{Guinaldo11}
M.~Guinaldo, D.~V. Dimarogonas, K.~H. Johansson, J.~Sanchez, and S.~Dormido,
  ``Distributed event-based control for interconnected linear systems,'' in
  \emph{50th IEEE Conference on Decision and Control and European Control
  Conference (CDC-ECC)}, 2011, pp. 2553--2558.

\bibitem{Cervin08CDC}
A.~Cervin and T.~Henningsson, ``Scheduling of event-triggered controllers on a
  shared network,'' in \emph{47th IEEE Conference on Decision and Control},
  2008, pp. 3601--3606.

\bibitem{Garcia12Parameter}
E.~Garcia and P.~J. Antsaklis, ``Parameter estimation in time-triggered and
  event-triggered model-based control of uncertain systems,''
  \emph{International Journal of Control}, vol.~85, no.~9, pp. 1327--1342,
  2012.

\bibitem{Guinaldo12CDC}
M.~Guinaldo, D.~Lehmann, J.~S{\'a}nchez, S.~Dormido, and K.~H. Johansson,
  ``Distributed event-triggered control with network delays and packet
  losses,'' in \emph{51st IEEE Conference on Decision and Control}, 2012, pp.
  1--6.

\bibitem{LiuHill12}
T.~Liu, D.~J. Hill, and B.~Liu, ``Synchronization of dynamical networks with
  distributed event-based communication,'' in \emph{51st IEEE Conference on
  Decision and Control}, 2012, pp. 7199--7204.

\bibitem{Garcia14CDCetof}
E.~Garcia and P.~J. Antsaklis, ``Event-triggered output feedback stabilization
  of networked systems with external disturbance,'' in \emph{53rd IEEE
  Conference on Decision and Control}, 2014, pp. 3566--3571.

\bibitem{WangLemmon08a}
X.~Wang and M.~Lemmon, ``Event triggered broadcasting across distributed
  networked control systems,'' in \emph{American Control Conference}, 2008, pp.
  3139--3144.

\bibitem{Mccourt14}
{M. J. McCourt, E. Garcia, and P. J. Antsaklis}, ``Model-based event-triggered
  control of nonlinear dissipative systems,'' in \emph{2014 American Control
  Conference}, 2014, pp. 5355--5360.

\bibitem{Adaldo15}
A.~Adaldo, D.~Liuzza, D.~V. Dimarogonas, and K.~H. Johansson, ``Control of
  multi-agent systems with event-triggered cloud access,'' in \emph{European
  Control Conference}, 2015, pp. 954--961.

\bibitem{nowzari19}
C.~Nowzari, E.~Garcia, and J.~Cortes, ``Event-triggered communication and
  control of networked systems for multi-agent consensus,'' \emph{Automatica},
  vol. 105, pp. 1--27, 2019.

\bibitem{Tallapragada12IFAC}
P.~Tallapragada and N.~Chopra, ``Event-triggered decentralized dynamic output
  feedback control for {LTI} systems,'' in \emph{IFAC Workshop on Distributed
  Estimation and Control in Networked Systems}, 2012, pp. 31--36.

\bibitem{Astrom02}
K.~J. Astrom and B.~M. Bernhardson, ``Comparison of riemann and lebesgue
  sampling for first order stochastic systems,'' in \emph{41st IEEE Conference
  on Decision and Control}, 2002, pp. 2011--2016.

\bibitem{Astrom08}
K.~J. Astrom, ``Event based control,'' in \emph{A. Astolfi and L. Marconi,
  (eds.), Analysis and Design of Nonlinear Control Systems}, Springer-Verlag,
  Berlin, 2008, pp. 127--147.

\bibitem{Wang06}
P.~Tabuada and X.~Wang, ``Preliminary results on state-triggered scheduling of
  stabilizing controls tasks,'' in \emph{45th IEEE Conference on Decision and
  Control}, 2006, pp. 282--287.

\bibitem{Tabuada07}
P.~Tabuada, ``Event-triggered real-time scheduling of stabilizing control
  tasks,'' \emph{IEEE Transactions on Automatic Control}, vol.~52, no.~9, pp.
  1680--1685, 2007.

\bibitem{Garcia13scl}
E.~Garcia and P.~J. Antsaklis, ``Output feedback networked control with
  persistent disturbance attenuation,'' \emph{Systems \& Control Letters},
  vol.~62, no.~10, pp. 943--948, 2013.

\bibitem{Garcia13MED}
{E. Garcia and P. J. Antsaklis}, ``Model-based control of continuous-time
  systems with limited intermittent feedback,'' in \emph{21st Mediterranean
  Conference on Control and Automation}, 2013, pp. 452--457.

\bibitem{Garcia12MED}
E.~Garcia and P.~J. Antsaklis, ``Model-based control of continuous-time and
  discrete-time systems with large network induced delays,'' in \emph{20th
  Mediterranean Conference on Control and Automation}, 2012, pp. 1135--1140.

\bibitem{GarciaBook}
E.~Garcia, P.~J. Antsaklis, and L.~A. Montestruque, \emph{Model-based control
  of networked systems}.\hskip 1em plus 0.5em minus 0.4em\relax Springer
  International Publishing, 2014.

\bibitem{lehmann2012event}
D.~Lehmann and J.~Lunze, ``Event-based control with communication delays and
  packet losses,'' \emph{International Journal of Control}, vol.~85, no.~5, pp.
  563--577, 2012.

\bibitem{Garcia20arxiv}
E.~Garcia and P.~J. Antsaklis, ``Model-based event-triggered control over lossy
  networks,'' \emph{arXiv paper no. 2007.15181}, 2020.

\end{thebibliography}

\end{document}